\documentclass{article}

\usepackage{fullpage}

\usepackage{graphicx}
\usepackage{graphics}
\usepackage{epsfig}
\usepackage{amsmath}
\usepackage{amssymb}
\usepackage[english]{babel}
\usepackage[latin1]{inputenc}
\usepackage{latexsym}

\newtheorem{proposition}{Proposition}
\newtheorem{lemma}{Lemma}
\newtheorem{remark}{Remark}

\def\R{\mathbb{R}}
\def\argmax{\mathop{\rm argmax}}
\def\ds{\displaystyle}
\def\qed{\hfill \ensuremath{\Box}}

\begin{document}

\title{Minimal-time bioremediation of natural water
resources}

\author{\sc P. Gajardo${}^{1}$, H. Ramírez${}^{2}$, A. Rapaport${}^{3}$\footnote{A. Rapaport and J. Harmand are with the Equipe-projet INRA-INRIA 'MODEMIC' (Modélisation et Optimisation des Dynamiques des Ecosystèmes MICrobiens)} and J. Harmand${}^{4}$\footnotemark[1]}
\date{\small
${}^{1}$ Departamento de Matemática, Universidad Técnica Federico Santa María. Avda. España 1680, Valparaíso, Chile.\\ Email: pedro.gajardo@usm.cl\\
${}^{2}$ Departamento de Ingeniería Matemática and Centro de Modelamiento Matemático (CNRS UMR 2807) Universidad de Chile. Avda Blanco Encalada 2120, Santiago,  Chile.\\ E-mail: hramirez@dim.uchile.cl\\
${}^{3}$ UMR 'MISTEA' Math\'ematiques, Informatique et STatistique pour l'Environnement et l'Agronomie (INRA/SupAgro) 2, place P. Viala, 34060 Montpellier,  France.\\ E-mail: rapaport@supagro.inra.fr\\
${}^{4}$ Laboratoire de Biotechnologie de l'Environnement, Route des Etangs, 11100 INRA Narbonne, France.\\ E-mail: harmand@supagro.inra.fr}

\maketitle

\begin{abstract} 
We study minimal time strategies for the treatment of pollution 
of large volumes, such as lakes or natural reservoirs, with the 
help of an autonomous bioreactor. The control consists in feeding the bioreactor from the resource, the clean output returning to the resource with the same flow rate. 
We first characterize the optimal policies among
constant and feedback controls, under the assumption of a uniform
concentration in the resource. 
In a second part, we study the influence of an inhomogeneity 
in the resource, considering two measurements points. With the help of the Maximum Principle, we show that the optimal control law is non-monotonic and terminates with a constant phase, contrary to the homogeneous case for which the optimal flow rate is decreasing with time. 
This study allows the decision makers to identify situations for 
which the benefit of using non-constant flow rates is significant.
\end{abstract}

\section{Introduction}

The fight against eutrophication of lakes and natural reservoirs 
(excessive development of phytoplankton associated with an excess of nutrients) constitutes a major challenge.
Such an ecological question has given rise to many studies over the last 30 years 
(see, for instance, the surveys \cite{gulati} or \cite{sondergaard} and references herein).
To remediate to eutrophication, many techniques such as bio-manipulation or ecological control have been proposed with mitigated results. 
A common point of the proposed remediation approaches is that they are usually based on "biotic" actions 
on the lake trophic chain dedicated to the restoration of the equilibrium of the local ecosystems. To do so, most studies are based on empirical knowledge. 
However, since the seventies, the use of eutrophication models (from heuristic data-based models at steady state to more recent dynamical mass-balance based models) together with optimal control techniques have been proposed (cf. \cite{estrada} and references herein).

\begin{center}
\begin{figure}[h]
\begin{center}
\includegraphics[width=7cm]{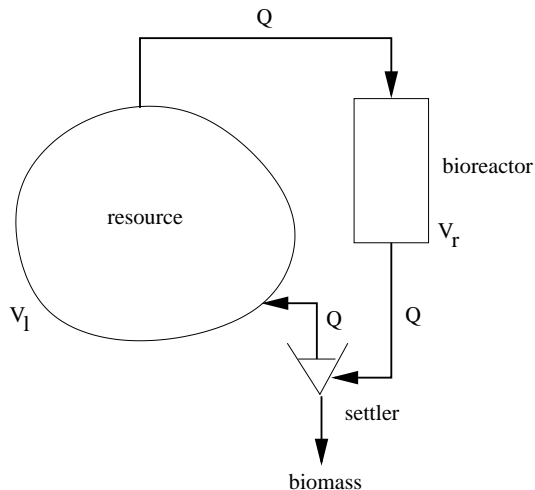}
\caption{Interconnection of the bioreactor with the resource.\label{Fig1}}
\end{center}
\end{figure}
\end{center}

In the present paper, an alternative to these techniques is studied using a very simple model 
of the lake. It is assumed that a small bioreactor is available to treat the polluted water
in removing a substrate considered as being in excess in the lake water.
More particularly, we consider a natural water resource of volume $V$
polluted with a substrate
of concentration $S_{l}$. As underlined above, typical examples of such natural water resources
to be treated are lakes or water tables that have been contaminated with
diffused pollutant as organic matter or nutrients.
The objective of the treatment is to make the concentration of such pollutant/contaminant $S_l$ decreasing down, as fast as possible,  to
a prescribed value $\underline S_{l}$, with the help of
a continuous stirred bioreactor of volume $V_r$.
The reactor is fed from the resource with a flow rate $Q$, and its output
returns to the resource with the same flow rate $Q$,
after separation of biomass and substrate in a settler
(see Figure \ref{Fig1}).
The settler avoids the presence of excessive biomass used for the treatment in the natural resource,
that could bring undesirable sludge and possibly lead to an increase of the eutrophication.
We assume that during the whole treatment, the volume
$V$ of the resource does not change.

Since the pioneer work by \cite{Ans},  the optimization of bioreactors operation
has received a great attention in the literature; see \cite{RaniRao,Banga,Alford} for reviews of the different optimization techniques that have been used in bioprocesses.
Among them, the theory of optimal control has proved to be a generic tool for deriving practical optimal rules \cite{VanImpeBastin,SmetsVanImpe,VanImpe}.
Clearly, one can distinguish two di\-ffe\-rent kinds of problems depending on the continuous
or discontinuous operation mode of the process. On one hand, if the process is operated in fed-batch, the
control objective is usually to optimize trajectories
for attaining a prescribed target in finite time or maximizing the production at a given time \cite{Hong,Lim,Johnson,Irvine,Kurtanjek,Shioya,Tsoneva,Moreno,Bonvin2003a,GHR08,Mazouni}. On the other hand, the optimal control of continuous processes usually involves a two steps procedure.
First, the optimal steady state is determined as a nominal set point,
maximizing a criterion \cite{Spitzer,Soukkou}. The benefit of ope\-ra\-ting a periodic control about the nominal point can be analyzed \cite{AbuleszLyberatos,RuanChen}.
Then, a control strategy that drives the state about the nominal
set point from any initial condition is searched for \cite{Kittisupakorn}, possibly in the presence of uncertainty on the model using extremum seeking techniques \cite{Wang,Zhang,Marcos,Mareels}

Concerning these strategies, the
problem studied in the present paper exhibits several original points  with respect to the contributions available in the literature. Indeed: 
\begin{itemize}
\item[-] The actual control problem is dedicated to the optimization
of transient trajectories - as in the case of fed-batch processes - while it is actually
a continuous process. It is due to the fact that in a standard optimal minimal-time problem of a bioprocess,
the volume of water to be processed is completely
decoupled from the bioreactor. In other terms,
the problem is to process, using a biological reactor, a given volume of "substrate" which is finally released
in the environment after processing (whatever it is operated continuously or discontinuously).
In the present problem,
the treated water is immediately recycled into the lake.
From the modelling point of view, this introduces an original coupling via the dilution of the treated water with the polluted one.
\item[-] The lake and the reactor are isolated in the sense no biomass is
supposed to be present in the water resource. The biomass used as a catalyst in the bioreactor
is separated from the treated water and withdrawn from the overall process.
Thus, in particular, the quantity of available biomass is not a limiting
parameter.
\end{itemize}

We consider the usual chemostat model for describing the dynamics 
of the bioreactor:
\begin{equation}
\label{eq1}
\left\{\begin{array}{lll}
\dot S_{r} & = & \ds -\mu(S_{r})X_{r} +\frac{Q}{V_r}(S_{l}-S_{r})\\
\dot X_{r} & = & \ds \mu(S_{r})X_{r} -\frac{Q}{V_r}X_{r}
\end{array}\right.
\end{equation}
where  $S_{r}$ and $X_{r}$ stand for the concentrations of
substrate and biomass, respectively.
For sake of simplicity, we assume that the yield coefficient of
this reaction is equal to one (at the price of changing the unitary
value of the biomass concentration, that is always possible).
The growth rate function $\mu(\cdot)$ fulfills the properties\\

\noindent {\bf Assumption A1.}\\
a. Function $\mu(\cdot)$ is increasing and such that $\mu(0)=0$.\\
b. Function $\mu(\cdot)$ is concave.\\

A reasonable hypothesis is to assume that the volume of the
resource is  much larger than the bioreactor one: $V >> V_r$,
and  that the possible variations of the manipulated variable 
$Q$ are slow 
compared to the time scale of the bioreactor dynamics. 
Consequently, one can consider that dynamics (\ref{eq1}) is
fast and its trajectories at the quasi-steady state 
$(S_{r}^{\star},X_{r}^{\star})=(S_{r}(Q),S_{l}-S_{r}(Q))$, where
$S_{r}(Q)$ fulfills $\mu(S_{r}(Q))=Q/V_r$ (see the usual equilibria
analysis of the chemostat \cite{SW95}).\\

\noindent {\bf Problem:}
The optimization problem consists in driving in minimal time the
concentration of the resource down to a prescribed value $\underline{S}_{l}>0$, playing with the control variable $Q>0$. In Section 2, we assume that this concentration is uniform in the resource, while in Section 3 we study the effect of a spatial inhomogeneity. 
For each case, we characterize the optimal policy
$Q^{\star}$ (resp. $Q^{opt}(\cdot)$) among constant (resp.
feedback controls). Section 4 is devoted to numerical simulations and discussions.

\section{The homogeneous case}

The dynamics of the resource concentration is simply
\begin{equation}
\label{reduced}
\dot S_{l} = \frac{Q}{V}(S_{r}(Q)-S_{l}).
\end{equation}
Notice that under Assumption A1.a, choosing
$Q$ is equivalent to choosing $S_{r}$ as a control variable:
\begin{equation}
\label{reducedS}
 \dot S_{l} = \alpha\mu(S_{r})(S_{r}-S_{l}), \quad
S_{r}\in (0,S_{l})
\end{equation}
where we denote  $\alpha=V_{r}/V$.

\begin{proposition}\label{prop1}
 Under Assumption A1, the best
cons\-tant control
$Q^{\star}$ is
defined as $Q^{\star}=V_r\mu(S_{r}^{\star})$, where $S_{r}^{\star}$ is the unique
minimum of the function
\begin{equation}
\label{Tf}
T_{f}(S_{r})=\frac{1}{\alpha\mu(S_{r})}\ln\left(\frac{S_{l}(0)-S_{r}}{\underline{S}_{l}-S_{r}}\right)
\end{equation}
on the interval $(0,\underline S_{l})$.
\end{proposition}

{\em Proof.} 
For a constant control, solutions of (\ref{reduced}) can be made
explicit:
\begin{equation}
S_{l}(t)= S_{r}(Q) + (S_{l}(0)-S_{r}(Q)) e^{-\frac{Q}{V}t},
\end{equation}
as well as the time $T_{f}(S_{r})$, given in (\ref{Tf}),
for reaching the target with $Q=V_{r}\mu(S_{r})$.
The function $T_{f}(\cdot)$ tends toward $+\infty$ when $S_{r}$ tends
toward $0$ or $\underline S_{l}$. Consequently, its infimum is reached
on the interval $(0,\underline S_{l})$.
Denote by $T_{f}^{\star}$ its minimum,
that we fix in the following. Then, for each constant control 
$S_{r}$, one has
\[
\begin{array}{ll}
\ds \frac{dS_{l}(T_f^{\star})}{d S_{r}} &=
1 - \left[1+\alpha\mu^{\prime}(S_{r})T_f^{\star}(S_{l}(0)-S_{r})\right]
e^{-\alpha\mu(S_{r})T_f^{\star}} \\[4mm]
\ds \frac{d^{2} S_{l}(T_f^{\star})}{d S_{r}^{2}}& =
\left[2\mu^{\prime}(S_{r})
+(\alpha\mu^{\prime}(S_{r})^{2}T_f^{\star} -\mu^{\prime\prime}(S_{r})(S_{l}(0)-S^{*}_r)\right]
\alpha T_f^{\star} e^{-\alpha\mu(S_{r})T_f^{\star}}
\end{array}
\]
and one deduces with Assumption A1 that the map
$S_{r} \mapsto S_{l}(T_{f}^{\star})$ is strictly convex.
Notice that one has necessarily $S_{l}(T_f^{\star})\geq \underline{S}_{l}$ and $S_{l}(T_f^{\star})= \underline{S}_{l}$ when
$S_{r}=S_{r}^{\star}$ realizes the  minimum of the
function $T_{f}(\cdot)$.
Consequently, the optimal control $S_{r}^{\star}$ is unique.\qed

\begin{proposition}\label{prop2} Under Assumption A1, the optimal feedback fulfills
$Q^{opt}(S_{l})=V_r\mu(S^{opt}_{r}(S_{l}))$
with
\begin{equation}
\label{feedback}
S^{opt}_{r}(S_{l}) \in
\argmax_{S_{r}\in (0,S_{l})}\mu(S_{r})(S_{l}-S_{r}) \ .
\end{equation}
Moreover, $t\mapsto Q^{opt}(t)$ is decreasing along any optimal trajectory.
\end{proposition}

{\em Proof.}
It is straightforward to check that the optimal feedback $S_{r}^{opt}$
is the one that makes the time derivative of $S_{l}$, given by
(\ref{reducedS}), the most negative at any time.
A necessary condition is to have
have $\mu^{\prime}(S^{opt}_{r})(S_{l}-S^{opt}_{r})=\mu(S^{opt}_{r})$.
Deriving this last expression w.r.t. time, one has
$\dot S^{opt}_{r}(2\mu^{\prime}(S^{opt}_{r})+\mu^{\prime\prime}(S^{opt}_{r})(S^{opt}_{r}-S_{l}))=\mu^{\prime}(S^{opt}_{r})\dot S_{l}$
and from Assumption A1, on obtains $\dot S^{opt}_{r}< 0$.\qed\\

For usual growth functions, one obtains the expressions:
\begin{center}
\begin{tabular}{c|c}
linear: $\mu(s)=\mu s$, &  Monod: $\ds \mu(s)=\frac{\mu_{\max}s}{K+s}$,\\[2mm]
\hline &\\
$\ds S^{opt}_{r}(S_{l})=S_{l}/2$ & 
$\ds S^{opt}_{r}(S_{l}) = \sqrt{K^{2}+KS_{l}}-K$
\end{tabular}
\end{center}

\section{Consideration of a spatial inhomogeneity}\label{2volumes}
\begin{center}
\begin{figure}[h]
\begin{center}
\includegraphics[width=7cm]{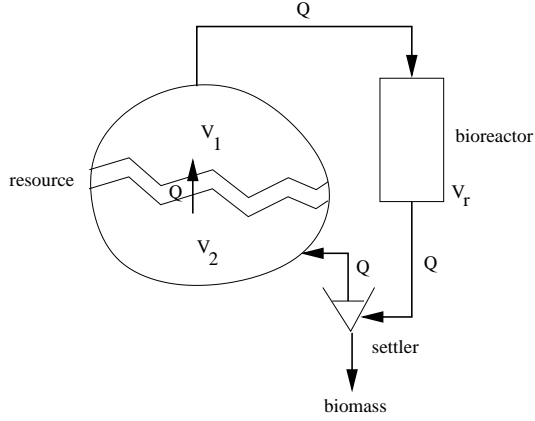}
\caption{Consideration of non-homogeneity in the resource.\label{Fig2}}
\end{center}
\end{figure}
\end{center}

The simplest way to introduce a gradient of
concentration in the model of the resource is to consider two compartments of volumes $V_{1}$, $V_{2}$ such that
$V=V_{1}+V_{2}$ (see Figure \ref{Fig2}), that we assume to be large with respect to $V_{r}$. Water is pumped from the first one while the clean one is rejected in the second one.
Denoting $\alpha_{i}=V_{r}/V_{i}$ ($i=1,2$), one obtains
the dynamics 
\begin{equation}
\label{reduced2}
\begin{array}{lll}
\dot S_{1} & = & \ds \frac{Q}{V_{1}}(S_{2}-S_{1})=
\alpha_{1}\mu(S_{r})(S_{2}-S_{1})\\[4mm]
\dot S_{2} & = & \ds \frac{Q}{V_{2}}(S_{r}(Q)-S_{2})=
\alpha_{2}\mu(S_{r})(S_{r}-S_{2})
\end{array}
\end{equation}
and can easily check that the domain
$\mathcal{D}= \{ (S_{1},S_{2})\in\R^{2}_{+} \, \vert \, S_{1}\geq
S_{2} \}$
is invariant for any control $S_{r}(\cdot)$ such that
$S_{r}(t) \in (0,S_{2}(t)]$ for any $t>0$.
We shall consider initial conditions in $\mathcal{D}$ and define the target 
$\mathcal{T} = \{ (S_{1},S_{2})\in \mathcal{D} \, \vert \, S_{1}\leq
\underline S_{l} \}$.

For $p\in[0,1]$ and $\tau\geq 0$, we define the function
\[
A(p,\tau)=\left|\begin{array}{ll}
\ds e^{-\alpha \tau} & \mbox{if } p=0\\[0mm]
\ds \frac{(1-p)e^{-\frac{\alpha\tau}{1-p}}
-p e^{-\frac{\alpha\tau}{p}}}{1-2p}
& \mbox{if } p\in (0,\frac{1}{2})\\[3mm]
\ds \left(1+\frac{\alpha}{2}\tau\right) e^{-\frac{\alpha\tau}{2}}
& \mbox{if } p=\frac{1}{2}\\[0mm]
A(1-p,\tau) & \mbox{if } p\in (\frac{1}{2},1]
\end{array}\right.
\]
and for $S_{0}>\underline{S}_{l}$ the function
\[
B(S_{r})=\frac{\underline{S}_{l}-S_{r}}{S_{0}-S_{r}} \ ,
\qquad S_{r} \in (0,\underline{S}_{l}) \ .
\]

\begin{proposition}
Let $p=V_{1}/V$.
For initial conditions such that $S_{1}(0)=S_{2}(0)=S_{0}>\underline S_{l}$, the best constant control $Q^{\star}$ and time $T_{f}^{\star}$
to reach the target ${\mathcal T}$ are defined by $Q^{\star}=V_{r}\mu(S_{r}^{\star})$, where $S_{r}^{\star}$ is such that the graph of $B(\cdot)$
touches tangentially the graph of
$S_{r} \mapsto A(p,\mu(S_{r})T_{f}^{\star})$ at
$S_{r}=S_{r}^{\star}$.
\end{proposition}

{\em Proof.}
The solution of (\ref{reduced2}) with constant
control $S_{r} \in(0,\underline S_{l})$ can be made explicit:
\[
S_{1}(t) = S_{r} + (S_{0}-S_{r})A(p,\mu(S_{r})t) \ ,
\]
and the time $T_{f}$ to reach the target fulfills
$A(p,\mu(S_{r})T_{f})= B(S_{r})$. Notice that one has the property
\[
S_{1}(t)\geq\underline{S}_{l} \quad \Longleftrightarrow \quad
A(p,\mu(S_{r})t)\geq B(S_{r}) \ ,
\]
from which the statement of the proposition follows.\qed\\

\begin{center}
\begin{figure}[h]
\begin{center}
\includegraphics[width=8cm]{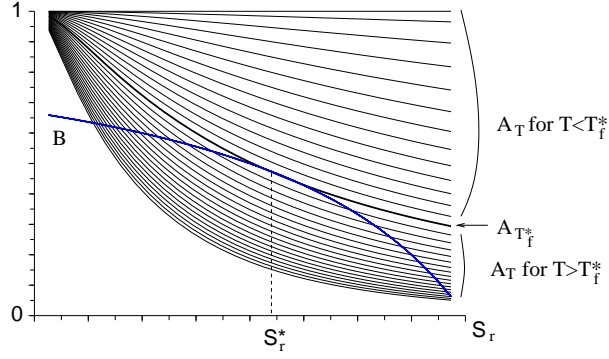}
\caption{Graphical determination of $S_{r}^{\star}$
and $T_{f}^{\star}$.\label{Fig3}}
\end{center}
\end{figure}
\end{center}

If one consider the family of functions
$A_{T}(S_{r})=A(p,\mu(S_{r})T)$,
parametrized by $T>0$, one has a gra\-phi\-cal interpretation 
of the optimum depicted in Figure  \ref{Fig3}.

$B(\cdot)$ is a concave function and one can check that 
$A_{T}(\cdot)$ are strictly
convex for values of $S_{r}$ large enough (under Assumption A1).
Consequently, there cannot exist more than one best constant control
 $S_{r}^{\star}$ in the domain where $A_{T_{f}^{\star}}$ is convex.

\begin{remark}
\label{remark1}
Notice that the case of an homogeneous resource
can be seen formally as the limiting case $p=0$ (although  cases when $p$ is close to $0$ are not compatible with assuming that $V_{1}$, $V_{2}$ are large with respect to $V_{r}$). 
One can easily check that $A(p,\tau) < A(0,\tau)$
for sufficiently large values of $\tau$.
Furthermore, times $T_{f}^{\star}$ are increasing with
respect to $S_{0}-\underline{S}_{l}$. 
Consequently, for initial conditions such that
$S_{0}$ is far from $\underline{S}_{l}$, the time
$T_{f}^{\star}$ is larger for an homogeneous resource
than for non-homogeneous one, even when the parameter $p$ is unknown
and the control $S_{r}^{\star}$ is determined for the
homogeneous case.
\end{remark}

For $S_1 > S_2> 0$ and  $\gamma> 0$, we define
\begin{eqnarray}
\ds \phi(S_1, S_{2},\gamma,S_r)&=&\mu(S_r)\left[1+\gamma\frac{(S_{2}-S_r)}{(S_1-S_2)}\right] \ , \label{def:phi}\\[3mm]
\ds \psi(S_1, S_{2},\gamma)&=&\mu^{\prime}(S_{2})-\gamma\frac{\mu(S_{2})}{S_1-S_2} \ . \label{def:psi}
\end{eqnarray}
The proof of the following lemma is left to the reader.
\begin{lemma}\label{lem3} Under Assumption A1, for $S_1 > S_2> 0$ and  $\gamma> 0$, the
function $\phi(S_1, S_{2},\gamma,\cdot)$ is strictly concave on $[0,S_{2}]$ and the property
\[
\max_{S_r\in(0,S_{2}]}\phi(S_1, S_{2},\gamma,S_r)=
\phi(S_1, S_{2},\gamma,S_{2}) \
\]
is fulfilled exactly when $\psi(S_1, S_{2},\gamma)\geq 0$.
\end{lemma}

\begin{proposition}\label{prop3}
 Under Assumption A1, from any initial
condition in $\mathcal{D}\setminus \mathcal{T}$, the optimal 
control $Q^{opt}(\cdot)$ consists in reaching a subset 
${\mathcal I}\subset\mathcal{D}\setminus \mathcal{T}$ from
which the constant control $Q^{opt}=V_{r}\mu(\underline S_{2})$ 
is optimal until $S_{1}(\cdot)$ reaches $\underline S_{l}$, where $\underline S_{2}$ is the value of $S_2$ when ${\mathcal I}$ is reached.
Moreover, $t \mapsto Q^{opt}(t)$ is increasing when approaching the set ${\mathcal I}$.
\end{proposition}

{\em Proof.}
Recall first that $\mathcal{D}$ is invariant.
If $S_{1}=S_{2}> \underline S_{l}$,
the feedback $S_r=S_{2}$  cannot be optimal (this would imply
$\dot S_{1}=\dot S_{2}=0$ at any time).
So,  any optimal trajectory is such that
$S_{1}(t)>S_{2}(t)$ for any $t>0$.\\
Let us write the Hamiltonian, along with the adjoint equations:
\[
H=-1+\max_{S_{r}\in[0,S_{2}]}
\mu(S_r)\left[\alpha_{1}\lambda_{1}(S_{2}-S_{1})+\alpha_{2}\lambda_{2}(S_r-S_{2})\right]
\]
\[
\left\{\begin{array}{ll}
\dot \lambda_{1} = \alpha_{1}\mu(S_r^{opt})\lambda_{1} \ , & \lambda_{1}(T_{opt})<0\\
\dot \lambda_{2} = \mu(S_r^{opt})(\alpha_{2}\lambda_{2}-\alpha_{1}\lambda_{1}) \ , \qquad & \lambda_{2}(T_{opt})=0
\end{array}\right.
\]
One deduce immediately that $\lambda_{1}(t)<0$ for any $t\geq 0$ and
can consider the function
\begin{equation}
\label{gamma}
\gamma(t)=\frac{\alpha_{2}\lambda_{2}(t)}{\alpha_{1}\lambda_{1}(t)}
\end{equation}
that fulfills 
$\dot \gamma = \mu(S_r^{opt})\left[(\alpha_{2}-\alpha_{1})\gamma -\alpha_{2}\right], \; \gamma(T_{opt})=0$.

Notice that  $\gamma=0$ implies $\dot \gamma <0$ and then
one obtains $\gamma(t)>0$ for any $ t \in [0,T_{opt})$.

When $S_{1}>S_{2}$, optimizing the Hamiltonian is equivalent to maximizing  $\phi(S_1, S_{2},\gamma,\cdot)$ (defined in \eqref{def:phi}), and then Lemma \ref{lem3} provides the uniqueness of $S_{r}^{opt}$.
A straightforward calculus gives
\begin{equation*}
 \frac{d}{dt}\left(\frac{\gamma}{S_{1}-S_{2}}\right)
=\frac{\alpha_2 \mu(S_r^{opt})}{S_1-S_2}\left[  \gamma \left(1-\frac{S_2-S_r^{opt}}{S_1 - S_2} \right)-1\right] .
\end{equation*}
From $\gamma(T_{opt})=0$ we deduce the existence of  $\tilde t \in [0,T_{opt})$ such that  $\frac{d}{dt}\left(\frac{\gamma}{S_{1}-S_{2}}\right) < 0$ for all $t \in [\tilde t, T_{opt}]$.
Then, for $\psi$ given by \eqref{def:psi}, one has
\begin{equation*}
\dot\psi=
\dot S_{2}\left(\mu^{\prime\prime}(S_{2})-\frac{\gamma}{S_1-S_2} \mu^{\prime}(S_{2})\right)-\frac{d}{dt}\left(\frac{\gamma}{S_{1}-S_{2}}\right) \mu(S_2),
\end{equation*}
which is positive
for $t \in [\tilde t, T_{opt}]$. Since  $\psi>0$
at $T_{opt}$, there exists $t_s \in [\tilde t,T_{opt})$ such that
$\psi\geq 0$ for $t\geq t_s$. Defining $\bar t=\inf\{t_s \in [\tilde t,T_{opt}) :\psi\geq 0 ~\mbox{for }t \in [t_s,T_{opt}]\}$,
with Lemma \ref{lem3}, one concludes that the optimal
$S_r^{opt}$ is constant equal to $\underline S_{2}=S_{2}(\bar t)$
at any time $t\geq \bar t$.\\
When $\bar t>0$, let us write $S_{r}^{opt}=uS_{2}$ with $u\in [0,1]$.
The left derivative $\dot u(\bar t^{-})$ has to be positive and
$\dot S_{2}(\bar t)=0$. This implies to have 
$\dot S_{r}^{opt}(\bar t^{-})>0$.\qed

\begin{proposition}
Under Assumption A1, for any initial
condition in $\mathcal{D}\setminus \mathcal{T}$, the optimal trajectory is unique.
\end{proposition}

{\em Proof.}
We recall, from the proof of Proposition \ref{prop3}, that 
along any optimal trajectory, one has
$S_{1}(t)>S_{2}(t)$ and $\gamma(t)>0$ for any $t \in(0,T_{f})$. 
Then, one has $\dot S_{1}<0$ and can re-parametrize the dynamics of 
variables $S_{2}$ and $\gamma$, defined in (\ref{gamma}), in terms of $S_{1}$ instead of time $t$, and write the non-autonomous dynamics
for optimal trajectories:
\begin{equation}
\label{reduced-gammaS1}
\begin{array}{llll}
\ds \frac{d S_{2}}{d S_1}& = & \ds \frac{ \alpha_2 }{\alpha_1} \frac{(S_{r}^{opt}-S_{2})}{(S_{2}-S_{1})}, &
\quad S_{2}(\underline S_{l})=S_{2}(T_{opt})\\[4mm]
\ds \frac{ d \gamma}{d S_1}  &=& \ds \frac{(\alpha_2-\alpha_1)\gamma-\alpha_2}{\alpha_1(S_{2}-S_{1})},
& \quad \gamma(\underline S_{l})=0
\end{array}
\end{equation}
where $S_{r}^{opt}$ is the unique maximum of $\phi(S_{1},S_{2},\gamma,\cdot)$ on $[0,S_{2}]$.
When $S_{r}^{opt}<S_{2}$ and $\gamma>0$, one has
\begin{equation}\label{eq:imply}
 \ds \frac{\partial \phi}{\partial S_r}=
\gamma\left(  \frac{S_2-S_{r}^{opt}}{S_1-S_2} -  \frac{ \mu(S_{r}^{opt})}{S_1-S_2} \right)+ \mu'(S_{r}^{opt}) =0\
\end{equation} and hence $ \ds \frac{\partial^{2} \phi}{\partial S_r\partial \gamma}=-\frac{\mu(S_{r}^{opt})}{\gamma}<0$.
\if{\begin{equation}\label{eq:imply}
\begin{split}
 \ds \frac{\partial \phi}{\partial S_r}=
\gamma\left(  \frac{S_2-S_{r}^{opt}}{S_1-S_2} -  \frac{ \mu(S_{r}^{opt})}{S_1-S_2} \right)+ \mu'(S_{r}^{opt}) =0\\
\Rightarrow  \ds \frac{\partial^{2} \phi}{\partial S_r\partial \gamma}=\frac{S_2-S_{r}^{opt}}{S_1-S_2} -  \frac{ \mu(S_{r}^{opt})}{S_1-S_2} =-\frac{\mu(S_{r}^{opt})}{\gamma}<0
\end{split}
\end{equation}}\fi
Fix $S_{1}$, $S_{2}$ and consider $S_{r}^{opt}$ as a function of $\gamma$. From \eqref{eq:imply}, one has
\[
\frac{\partial^{2} \phi}{\partial S_r\partial \gamma}+\frac{\partial^{2} \phi}{\partial S_r^{2}}\frac{\partial S_{r}^{opt}}{\partial \gamma}=0
\]
and from the strict concavity of $\phi(S_{1},S_{2},\gamma,\cdot)$, 
given by Lemma \ref{lem3}, one deduces  $\frac{\partial S_r^{opt}}{\partial \gamma} \le 0$.

The Jacobian matrix of system \eqref{reduced-gammaS1} is of the form
\[
\left[\begin{array}{cc}
\star & \frac{ \alpha_2 }{\alpha_1(S_{2}-S_{1})} \frac{\partial S_{r}^{opt}}{\partial \gamma}\\[3mm]
-\frac{((\alpha_2-\alpha_1)\gamma-\alpha_2)}{\alpha_1(S_{2}-S_{1})^2} & \star
\end{array}\right]
\]
from which one observes the non-negativity of off-diagonal terms, because $(\alpha_2-\alpha_1)\gamma-\alpha_2=\dot \gamma/\mu(S_r^{opt}) < 0$.
So, the dynamics (\ref{reduced-gammaS1}) is cooperative (in time $S_{1}$), and since $\gamma(\underline{S}_{l})=0$, one deduces that 
two solutions of (\ref{reduced-gammaS1}) cannot cross in the $(S_{1},S_{2})$ plane. Finally, one obtains the uniqueness of the optimal trajectory for a given initial condition in 
$\mathcal{D}\setminus \mathcal{T}$.\qed

\begin{remark}
\label{remark2}
When $S_{1}(0)=S_{2}(0)=S_{l}$, from the expression of the Hamiltonian we obtain that the 
optimal control $S_r^{opt}$ is such that at the beginning  it maximizes
$\mu(S_{r})(S_{l}-S_{r})$. Therefore, it is exactly the same as in the 
homogeneous case of Proposition \ref{prop2}.
Measuring the initial rate of variation of $S_{2}$
gives an estimation of the parameter $\alpha_{2}$ to
fit the model, as one has
\[
\dot S_{2}(0)=\alpha_{2}\mu(S_{r}^{opt})(S_{2}(0)-S_{r}^{opt}) \ .
\]
If it is close to $\alpha$, then the model with one compartment should suit.
\end{remark}

For $S_{1}^{0}>\underline S_{l}>S_{2}^{0}$, we define
when $\alpha_{1}\neq \alpha_{2}$:
\[
f_{0}(S_{1}^{0},S_{2}^{0}) =\frac{ \alpha_2\left(1-\left(\frac{S_{1}^0-S_{2}^0}{\underline
  S_{l}-S_{2}^0}\right){\frac{\alpha_1-\alpha_2}{\alpha_1}} \right)}{(\alpha_2-\alpha_1)(S_1^0- S_2^0)}, \;
\beta=\left(\frac{\alpha_1}{\alpha_2} \right)^{\frac{\alpha_1}{\alpha_1-\alpha_2}}
\]
and when $\alpha_{1}=\alpha_{2}$:
$\ds f_{0}(S_{1}^{0},S_{2}^{0}) =\frac{\ln\left(\frac{S_{1}^0-S_{2}^0}{\underline
  S_{l}-S_{2}^0}\right)}{(S_1^0- S_2^0)}, \; \beta=e$.

\begin{proposition}
\label{prop-final-stage}
The set ${\mathcal I}$, where a constant control is optimal,  is given by
\[
{\mathcal I} = \{ (S_{1}^{0},S_{2}^{0})\in (\underline S_{l},+\infty)\times (0,\underline{S}_{l}) \mbox{ s.t. } S_{2}^{0} \leq \bar S_{2} ~\mbox{ or } \mu(S_{2}^{0})f_{0}(S_{1}^{0},S_{2}^{0}) \leq \mu^{\prime}(S_{2}^{0})\,\}
\]
where $\bar S_{2}$ is the unique solution in $(0,\underline{S}_{l})$ of
\begin{equation}
\label{intersection}
\mu(\bar S_{2})=\beta\mu^{\prime}(\bar S_{2})(\underline{S}_{l}-\bar S_{2}) \ .
\end{equation}
\end{proposition}

{\em Proof.}
With control $Q=V_{r}\mu(S_{2}^{0})$, $S_{2}(\cdot)$ is equal to $S_{2}^{0}$ and solution $S_{1}(\cdot)$ can be made explicit. Then, time $T_{f}$ such that $S_{1}(T_{f})=\underline{S}_{l}$, and solution $\gamma(\cdot)$ such that $\gamma(T_{f})=0 $ can be also made explicit. Let $f(t)=\gamma(t)/(S_{1}(t)-S_{2}^{0})$. According to Lemma \ref{lem3}, this constant strategy is optimal exactly when
\begin{equation}
\label{condC}
\mu^{\prime}(S_{2}^{0})\geq \mu(S_{2}^{0})f(t) \ , \quad
t \in [0,T_{f}] \ . 
\end{equation}
One can easily check that $\dot f =\alpha_{2}\mu(S_{2}^{0})(\gamma-1)/(S_{1}-S_{2}^{0})$, and consequently, $\dot f$ cannot be null more than one time (recall from the proof of Proposition \ref{prop3} that $\gamma$ 
is non-increasing). One has also $f(0)>0$, $f(T_{f})=0$, and $f^{\prime}(T_{f})\leq 0$.

\begin{center}
\begin{figure}[h]
\begin{center}
\includegraphics[width=8cm]{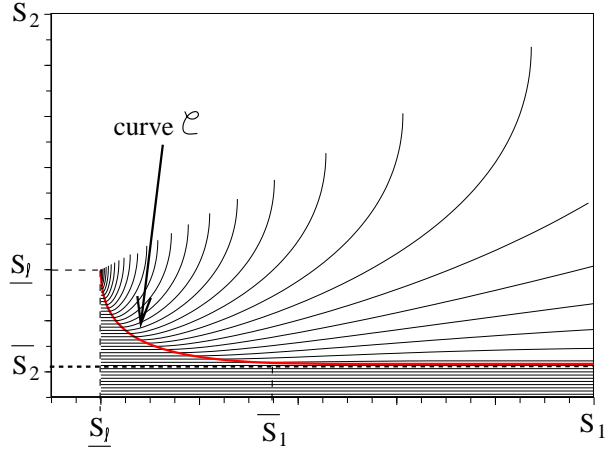}
\caption{Backward integration of the extremals.\label{Fig4}}
\end{center}
\end{figure}
\end{center}

If $f^{\prime}(0)\leq 0$, then condition (\ref{condC}) is equivalent
to have $(S_{1}^{0},S_{2}^{0})$ below the graph of the curve ${\mathcal C}$ defined by $\mu(S_{2}^{0})f_{0}(S_{1}^{0},S_{2}^{0})= \mu^{\prime}(S_{2}^{0})$. A straightforward but lengthy computation gives
$f(0)=f_{0}(S_{1}^{0},S_{2}^{0})$). One can also check that 
$f^{\prime}(0)\leq 0$ is equivalent to have $(S_{1}^{0},S_{2}^{0})$ below the line ${\mathcal L}$ defined by $S_{1}^{0}-S_{2}^{0}=\beta(\underline{S}_{l}-S_{2}^{0})$.
The intersection point $(\bar S_{1},\bar S_{2})$ of ${\mathcal C}$ and ${\mathcal L}$ is given by $\bar S_{2}$ solution of (\ref{intersection}), its uniqueness being guaranteed by the concavity of $\mu$. One can easily check that ${\mathcal C}$ is below ${\mathcal L}$ for 
any $S_{2}\in [\bar S_{2},\underline{S}_{l}]$.\\
If $f^{\prime}(0)\geq 0$, on can check that 
$\max_{t}f(t)=1/\beta(\underline{S}_{l}-S_{2}^{0})$ and then 
condition (\ref{condC}) is equivalent to have 
$S_{2}^{0}\leq \bar S_{2}$. 
Moreover, the straight  line $S_2=\overline S_2$ is below  
the graph of the curve $\mathcal{C}$ in the interval 
$[\underline S_l, \overline S_1]$ (see Figure \ref{Fig4}).\qed

For the Monod law, one can find
\[
\overline S_2= \frac{1}{2}\left( -K(1+\beta) + \sqrt{K^2(1+\beta)^2 + 4K\beta \underline S_l}\right) \ .
\]

\section{Discussion and numerical simulations}

The benefits of our theoretical analysis
is to identify efficient pumping strategies and some 
of their robustness properties. 
We summarize those contributions in terms of
the following rules for the decision makers:\\
1. The profit of using the optimal feedback strategy compared 
to the best constant one, can be easily determined numerically
(see the simulations below). As expected, the more the resource is initially polluted, the more the 
improvement of the feedback policy is significant. 
Depending on the ratio "initial pollution over desired maximal pollution level'', the decision maker can then decide 
whether it worth adopting a time-varying strategy.\\
2. A spatial inhomogeneity of the pollution concentration improves the treatment time on the condition that the resource is enough polluted. Moreover, applying the best constant strategy as if the resource was perfectly homogeneous is robust with respect to uncertainty on the inhomogeneity parameter in the sense that it provides a guaranteed time (see Remark \ref{remark1}).\\
3. Measuring the initial speed of variation of concentration at two
remote locations in the resource allows to identify the inhomogeneity parameter of the model (see Remark \ref{remark2}). Then, the decision maker can decide if it worth considering a feedback strategy with two measurement points instead of one.\\
4. The optimal feedback strategy for the inhomogeneous case
consists in applying a constant flow rate when 
the concentration $S_{2}$ reaches a prescribed value given by 
Proposition \ref{prop-final-stage}. The concentration 
$S_{2}$ is then maintained constant, without having to measure 
the concentration $S_{1}$ (see Figure \ref{Fig4}).\\

\begin{center}
\begin{figure}[h]
\begin{center}
\includegraphics[width=8cm]{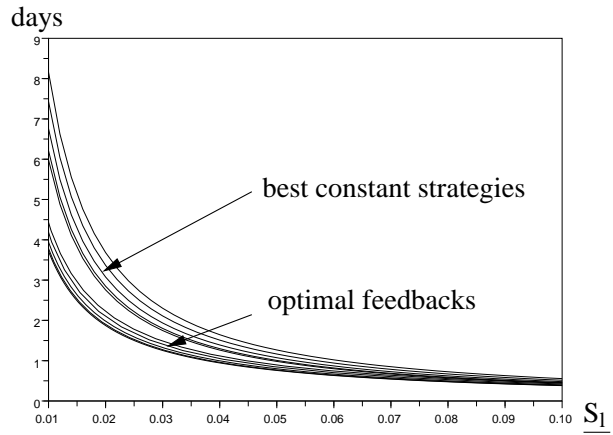}
\caption{Comparison of optimal times.\label{Fig5}}
\end{center}
\end{figure}
\end{center}

\begin{center}
\begin{figure}[h]
\begin{center}
\includegraphics[width=8cm]{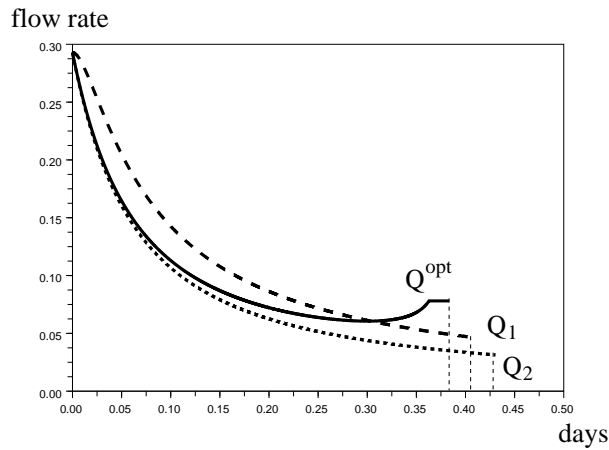}
\caption{Optimal and sub-optimal policies.\label{Fig6}}
\end{center}
\end{figure}
\end{center}

Simulations have been conducted for the Monod law with $\mu_{\max}=1 \,s^{-1}$, $K=1 \, mol.m^{-3}$, and volumes $V=1000\, m^{3}$,  $V_{r}=1\, m^{3}$. The initial concentration of pollutant has been chosen equal to $1 \, mol.m^{-3}$ uniformly. Figure \ref{Fig5} shows the comparison of minimal times for different values of $\underline S_{l}$ (the curves corresponds to different values of the parameter $p$).
On this example, one can see that for $\underline S_{l}=0.01$, 
the mi\-ni\-mal time among constant controls is about twice larger 
than among feedbacks. The influence of inhomogeneity 
is also quite significant.\\
The optimal feedback (\ref{feedback}) for the homogeneous case is a simple law that 
provides a decreasing flow rate $Q$ w.r.t. time (cf. Proposition \ref{prop2}), contrary to the inhomogeneous case for which it is
non-monotonic (cf. Proposition \ref{prop3}). 
In Figure \ref{Fig6}, we compare the optimal policy $Q^{opt}(\cdot)$ 
for $p=0.4$ and $\underline S_{l}=0.1$ with 
$Q_{1}(\cdot)$, resp. $Q_{2}(\cdot)$  applying 
the formula (\ref{feedback}) on measurement $S_{1}$, resp. $S_{2}$.

The true optimal feedback control, in the model  that consider both measurements, is more sophisticated in the sense that 
it anticipates the approach to the target, increasing the flow rate and freezing it.\\  
The study has been made assuming that the steady-state characteristics $Q \mapsto S_{r}(Q)$ of the bioreactor is perfectly known. 
Uncertainty on this map as well as on measurements will be the matter of a forthcoming work.

\section*{Acknowledgment}
This research was developed in the context of DYMECOS INRIA associated team, being partially supported by INRIA-CONICYT French-Chilean cooperation program and MOMARE SticAmsud project .
The first and third author thank the support of FONDECYT projects 1080173 and
1070297, and
Fondo Basal, Centro de Modelamiento Matematico - U.
de Chile. 

\bibliographystyle{plain} 
\bibliography{biblio_lake}
\end{document}